\documentclass[12pt,reqno]{amsart}

\usepackage{setspace}
\usepackage{apalike}
\usepackage[arrow,matrix,curve]{xy}

\usepackage{amssymb, latexsym, amsmath, amscd, array, hyperref
}

\usepackage{breakurl}   

\usepackage{graphicx} 

\usepackage[export]{adjustbox}  

\usepackage{epigraph}

\theoremstyle{definition}

\numberwithin{equation}{section}

\DeclareMathOperator{\adequal}{\,{}_{\ulcorner\!\urcorner}\,}

\newcommand\R {{\mathbb R}}

\author{Mikhail G. Katz} \address{Department of Mathematics, Bar Ilan
  University, Ramat Gan 5290002 Israel} \email{katzmik@math.biu.ac.il}

\author{Karl Kuhlemann} \address{Gottfried Wilhelm Leibniz University
  Hannover, D-30167 Hannover, Germany}
\email{kus.kuhlemann@t-online.de}

\author{Semen S. Kutateladze} \address{Sobolev Institute of
  Mathematics, Novosibirsk State University, Russia}

\begin{document}

\thispagestyle{empty}


\title [Marx versus Engels on infinitesimals] {Marx versus Engels on
  infinitesimals: Chimera or triumph?}



\begin{abstract}
We document the evolution of Karl Marx's take on infinitesimals.  We
contrast his initial favorable stance with later criticisms, and
examine the differing perspectives of Marx and Engels on the subject.
Marx's favorable assessment was based on his study of Sauri's
textbook.  Later, influenced by Boucharlat's textbook, Marx reversed
his position to an unfavorable stance, describing belief in
infinitesimals as a `chimera'.  Marxist scholar Guglielmo Carchedi
claims that ``Marx differentiates with the eyes of the social
scientist, of the dialectician'' but fails to note dialectician
Engels' endorsement of infinitesimals.

Struik linked Marx to Abraham Robinson, but missed the fact that the
link passes via{\ldots}~Fermat.  Namely, there may be an affinity, as
per Struik, between Marx's comments on the calculus and Robinson's
nonstandard analysis, but the kernel of such an affinity resides in
the techniques already found in the context of Fermat's
\emph{adequality}.  To adapt Carchedi's metaphor, we could say that
Marx may have differentiated with the eyes of \emph{adaequo} of Pierre
de Fermat.

The first editor who worked on some of Marx's mathematical manuscripts
in the mid-1920s was Emil J. Gumbel, though he is not mentioned in
either the 1933 or the 1968 Soviet edition of Marx's mathematical
manuscripts.
\end{abstract}

\maketitle


\epigraph{Marx differentiates with the eyes of the social scientist,
  of the dialectician.  --Carchedi}

\epigraph{Your remark that: \emph{``However (this comparison of the
  Taylor's Series with the Binomial Theorem) was not put forward by
[Marx], it is rather customary''} is quite unnecessary.  The Institute
  of Marx and Engels thinks that it would not be expedient to publish
  such remarks in print, so the best course for us would be to strike
  out this gloss.  \\ --Alexei Voden to Emil J. Gumbel}

\section{Differentiation, dialectics, and adequality}
\label{s1}

There has been a flurry of recent publications on the mathematical
manuscripts of Karl Marx, including
\cite{Ri18},
\cite{Ki19},
\cite{Da20},
\cite{Br21},
\cite{Ma21},
\cite{Ol21},
and
\cite{Vo24}.%
\footnote{Our second epigraph is from
\cite[31]{Vo19}.}
We propose to take a fresh look at the matter.

Marx seems to have had the distinction of being one of the first to
have applied the term \emph{chimera} (in the sense of an
\emph{illusory thing}) to disparage belief in infinitesimals.
Controversy in this area existed at least since the inception of the
infinitesimal calculus in the seventeenth century, but the specific
term \emph{chimera} was first used to denigrate infinitesimals by
d'Alembert and Moigno;%
\footnote{\label{f2}According to \cite[75]{Ma83}, Marx was in
possession of Moigno's work \cite{Mo40} where the term was used in
reference to infinitesimals.}
see \cite{25a}.
Marx commented as follows: ``The closely held belief of some
rationalising mathematicians that~$dy$ and~$dx$ are quantitatively
actually only infinitely small, only approaching 0/0, is a chimera''
\cite[5]{Ma83}; see also \cite[260]{Fa09}.
In a similar vein, Marx described the infinitesimal calculus of
Leibniz and Newton as \emph{mystical} (see Section~\ref{s2}).  In this
Marx was followed by a number of authors, not all of whom may have
been aware of the source of the modifier `mystical'.  Thus, Richard
Courant noted that ``Even in the case of scientists of the first rank
we find operations with the new concepts based chiefly on a feeling
for the right result and not always free from \emph{mystical}
associations -- particularly in the case of the ominous `infinitely
small quantities' or \mbox{`infinitesimals'}'' \cite[56; emphasis
  added]{Cou37}.  Moreover, ``in the early days of the differential
calculus even Leibnitz%
\footnote{\label{cou}Courant's spelling.}
himself was capable of combining these vague \emph{mystical} ideas
with a thoroughly clear understanding of the limiting process.  It is
true that this fog which hung round the foundations of the new science
did not prevent Leibnitz or his great successors from finding the
right path.  But this does not release us from the duty of avoiding
every such hazy idea in our building-up of the differential and
integral calculus'' \cite[81; emphasis added]{Cou37}.

\cite{Mc} explores Courant's engagement with socialism and marxism.
As late as 1952, shortly after being elected to the Hungarian Academy
of Sciences, the logician L\'aszl\'o Kalm\'ar ``exhorts Hungarian
physicists not to use capitalist infinitesimals in their work but to
employ the socialistic~$\epsilon$-$\delta$ method in their proofs'';
see \cite{Ke59}, \cite{Ma91}.  Carchedi for his part queries: ``Which
view of social reality is hidden behind and informs Marx's method of
differentiation?  Marx differentiates with the eyes of the social
scientist, of the dialectician'' \cite[424]{Car}, \cite[289]{Ca11}.

Carchedi argues that Marx differentiated with the eyes of a
dialectician.  Did he really?  We argue that Karl Marx (1818--1883)
differentiated with the eyes of the \emph{adaequo} (adequality) of
Pierre de Fermat (1601/7--1665).  Marx's excursion into higher
mathematics is reminiscent of the procedures introduced by the
seventeenth century master (several decades before the calculus was
officially propounded by Leibniz and Newton), as we elaborate
in Section~\ref{s6}.

In his review of
\cite{Ke77}, Dirk Jan Struik linked Marx to Abraham Robinson: ``Marx's
methods invite comparison with the present non-standard analysis of
A. Robinson or the reintroduction of infinitesimals into the
classroom''\cite{St77}.  Yet Struik, following \cite[316]{Ke77},
failed to appreciate the fact that the link passes via seventeenth
century work.  Namely, there is indeed an affinity between Marx's
comments on the calculus and Robinson's nonstandard analysis.
However, the kernel of such an affinity resides in the techniques
already found in the context of Fermat's adequality.

In Sections \ref{s2} through \ref{s4}, we compare the attitudes toward
infinitesimals in Marx and Engels.  In Section~\ref{s3}, we analyze
Marx's position in more detail.  In Section~\ref{s7}, we explore
Kolman's role in the hagiography of Marx's mathematics.
Section~\ref{s8} details the history of the attempts to edit Marx's
mathematical manuscripts, by Gumbel, Yanovskaya, Antonova, and others.
Section~\ref{s9} reviews the tumultuous history of Leibnizian
infinitesimals.

\section{Thesis: Mystical Differential Calculus}
\label{s2}

Ever the systematizer, Marx distinguished three periods in the
development of the calculus:
\begin{enumerate}
\item
  the `mystical differential calculus' of Newton and Leibnitz,%
\footnote{Marx's spelling.}
\item
the `rational differential calculus' of Euler and d'Alembert,
\item
the `pure algebraic calculus' of Lagrange;
\end{enumerate}
see \cite[xx]{Ma83}.\, Marx's comments on the first period were
characteristically cutting: ``Thus: they themselves believed in the
mysterious character of the newly discovered calculus, that yielded
true (and moreover, particularly in the geometrical application,
astonishing) results by a \emph{positively false} mathematical
procedure. They were thus self-mystified, valued the new discovery all
the higher, enraged the crowd of old orthodox mathematicians all the
more, and thus called forth the cry of opposition, that even in the
lay world has an echo and is necessary in order to pave the way for
something new'' \cite[168; emphasis added]{Ma68}.  However, the
procedure in question was not `positively false'.  The consensus among
historians today is that the Leibnizian calculus was far more coherent
than implied by the above passage, though there is disagreement as to
the exact nature of Leibniz's fictional infinitesimals.  For an
analysis, see \cite{23f}.

In Section~\ref{two}, we will see that Engels held a more charitable
view of the calculus of Newton and Leibniz.  Andrea Ricci correctly
notes that ``[Marx's] aim was to derive the derivative directly from
the process of variation of the function so that its algebraic
{\ldots}~origin is met'' \cite[Abstract]{Ri18}.  Such an algebraic
approach works for polynomials and possibly analytic functions, but
not necessarily for more general functions.

\section{Antithesis: Engels on infinitesimals as a triumph}
\label{two}

Friedrich Engels felt rather differently from Karl Marx about
infinitesimals.  He viewed infinitesimals as foundational in
understanding nature, contrasting with Marx's scepticism.

\subsection{\emph{Dialectics of Nature} on infinitesimals}

In his manuscript \emph{Dialectics of Nature}, Engels wrote that~$dx$
is infinitely small, but effective and generates everything, or in
more detail: ``Identity and difference -- the dialectical relationship
that already appears in differential calculus, where~$dx$ is
infinitely small but still effective and generates everything.''%
\footnote{``Identit\"at und Unterschied -- das dialektische
Verh\"altnis schon in der Differentialrechnung, wo~$dx$ unendlich
klein, aber doch wirksam und alles macht'' 
\cite{En}.}
The manuscript of \emph{Dialectics of Nature} includes passages in
Marx's handwriting with quotations from Greek philosophers, but the
passages on infinitesimals are in Engels' own handwriting.  Engels
writes further:
\begin{quote}
Let us take an example.  Of all theoretical advances, probably none is
regarded as such an exalted triumph of the human spirit as the
invention of infinitesimal calculus in the last half of the 17th
%
%
century.  If anywhere, we have here a pure and exclusive act of the
human spirit.  The mystery that still surrounds the quantities used in
the infinitesimal calculus -- the differentials and infinities of
various degrees -- is the best proof that we still imagine that we are
dealing here with pure “free creations and imaginations” of the human
spirit, for which the objective world offers no equivalent.%
\footnote{``Nehmen wir ein Beispiel.  Von allen theoretischen
Fortschritten gilt wohl keiner als ein so hoher Triumph des
menschlichen Geistes wie die Erfindung der Infinitesimalrechnung in
der letzten H\"alfte des 17. Jahrhunderts.  Wenn irgendwo, so haben
wir hier eine reine und ausschlie{\ss}liche Tat des menschlichen
Geistes.  Das Mysterium, das die bei der Infinitesimalrechnung
angewandten Gr\"o{\ss}en -- die Differentiale und Unendlichen
verschiedener Grade -- noch heute umgibt, ist der beste Beweis
daf\"ur, da{\ss} man sich noch immer einbildet, man habe es hier mit
reinen `freien Sch\"opfungen und Imaginationen' des Menschengeistes zu
tun, wof\"ur die objektive Welt kein Entsprechendes biete.
\cite{Enzeno}; cf.~\cite[218]{En}.}
\end{quote}
His conclusion is unequivocal:
``And yet the opposite is the case. Nature provides the models for all
these imaginary quantities.''%
\footnote{``Und doch ist das Gegenteil der Fall.  F\"ur alle diese
imagin\"aren Gr\"o{\ss}en bietet die Natur die Vorbilder.'' (ibid.)}
Engels' respect for infinitesimals may have been inspired by a
contemporary neo-Kantian philosopher Hermann Cohen (despite their
philosophical differences in the context of the materialist/idealist
divide; see further in Section~\ref{s71b}).  \cite{Co83} spoke of
infinitesimals in similar terms (though he would have balked at
`dialectical' accounts for infinitesimals) and placed them at the root
of his philosophical system for the natural sciences.  On Cohen and
his students, see \cite{Mo03}, \cite{13h}, \cite{Ed20}, \cite{Pr23}.
Giovanelli documents Cohen's defense of infinitesimals, posited as an
alternative to limits, as follows:
\begin{quote}
Cohen seems to include in the category of the `method of limits' all
attempts to justify the calculus without using any methods extraneous
to traditional algebra {\ldots}~For Cohen, these approaches represent
a sort of `repression' of a concept that was apparently unbearable to
the mathematical consciousness, the concept of the `infinitesimally
small,'~{\ldots}~he claims that such justifications conceal the
authentic motive that induced the `discoverers' of the infinitesimal
`method' to introduce a new type of magnitudes unknown to previous
mathematics: `This reading obscures the discovery and its tendency;
the positive, autonomous, irreplaceable element that forms the basis
of this new type of magnitudes is smoothed out [nivellirt]'.
\cite[11]{Gi16}.
\end{quote}
Giovanelli's quotation of Cohen is translated from \cite[95]{Co83}.
Cohen had been developing the ideas of the book since 1881; see
\cite[296]{Gi17}.

\subsection{\emph{Anti-D\"uhring} on infinitesimals}

Engels defended infinitesimals not only in \emph{Dialectics of Nature}
but also in his published work known as \emph{Anti-D\"uhring}, where
he wrote: ``Elementary mathematics, the mathematics of constant
quantities, moves within the confines of formal logic, at any rate on
the whole; the mathematics of variables, whose most important part is
the infinitesimal calculus, is in essence nothing other than the
application of \emph{dialectics} to mathematical relations''
\cite[emphasis added]{En78}.  Engels went on to mention Leibniz's
efforts to defend the infinitesimal calculus in his time:
\begin{quote}
To attempt to prove anything by means of \emph{dialectics} alone to a
crass metaphysician like Herr D\"uhring would be as much a waste of
time as was the attempt made by Leibniz and his pupils to prove the
principles of the infinitesimal calculus to the mathematicians of
their time.  The differential gave them the same cramps as Herr
D\"uhring gets from the \emph{negation of the negation}, in which,
moreover, as we shall see, the differential also plays a certain role
(Engels, ibid.).
\end{quote}
We will refrain from commenting on Engels' \emph{ad hominem} attacks
above, other than noting that neither D\"uhring nor Engels possessed a
viable approach to the infinitesimal calculus, the difference being
that D\"uhring never mendaciously claimed otherwise; for further
details see \cite{He85}.

What is notable here is Engels' emphasis on dialectics and the concept
of the negation of the negation in ``proving the principles'' of the
infinitesimal calculus.  On the rhetoric of \emph{dialectics} and
\emph{negation of the negation}, see further in Section~\ref{s71}.

\subsection{Analysis of Engels' remarks}

It is significant that Engels spoke in an approving way specifically
about differentials and infinities (for which ``Nature provides
models''), rather than merely about the (infinitesimal) calculus as a
mathematical technique.  To Engels, they clearly were not to be
dismissed as \emph{chimeras}.

The fact that Engels stresses the seventeenth century origin of the
infinitesimal calculus is further evidence that what he had in mind
was the infinitesimal calculus as conceived by its inventors, namely
using genuine infinitesimals (as opposed to nineteenth century efforts
to reformulate the calculus without them).

Similarly significant is his postulation of infinitesimals and
infinities of various orders in Nature (``Nature provides the models
for all these imaginary quantities'').
Engels provides the following clarification:

\begin{quote}
Our geometry is based on spatial relationships, our arithmetic and
algebra on numerical quantities that correspond to our earthly
relationships, which therefore correspond to the body sizes that
mechanics calls masses -- masses as they occur on Earth and are moved
by people.  Compared to these masses, the mass of the Earth appears
infinitely large and is also treated as infinitely large by earthly
mechanics. Earth's radius~$= \infty$, the principle of all mechanics
in the law of gravity.  However, not only the Earth, but also the
entire solar system and the distances within it appear to be
infinitely small as soon as we consider the distances in the star
systems visible to us telescopically, which can be estimated in light
years.  So here we already have an infinite not only of the first but
also of the second degree, and we can leave it to the imagination of
our readers to construct further infinities of higher degrees in
infinite space, if they feel like it.%
\footnote{``Unsre Geometrie geht aus von Raumverh\"altnissen, unsre
Arithmetik und Algebra von Zahlengr\"o{\ss}en, die unsren irdischen
Verh\"altnissen entsprechen, die also den K\"orpergr\"o{\ss}en
entsprechen, die die Mechanik Massen nennt – Massen, wie sie auf der
Erde vor kommen und von Menschen bewegt werden.  Gegen\"uber diesen
Massen erscheint die Masse der Erde unendlich gro{\ss} und wird von
der irdischen Mechanik auch als unendlich gro{\ss} behandelt.
Erdradius~$= \infty$, Grundsatz aller Mechanik im Fallgesetz. Aber
nicht nur die Erde, sondern das ganze Sonnensystem und die in ihm
vorkommenden Entfernungen erscheinen ihrerseits wieder als unendlich
klein, sobald wir uns mit den nach Lichtjahren zu schätzenden
Entfernungen in dem für uns teleskopisch sichtbaren Sternensystem
besch\"aftigen.  Wir haben hier also schon ein Unendliches nicht nur
des ersten, sondern auch des zweiten Grades, und k\"onnen es der
Phantasie unsrer Leser \"uberlassen, sich noch weitere Unendliche
h\"oherer Grade im unendlichen Raum zurechtzukonstruieren, falls sie
dazu Lust versp\"uren.'' 
\cite{En}.}
\end{quote}
There are interesting similarities between this passage from Engels
and a passage from Leibniz's 1702 letter where he sought to motivate
the use of infinitesimals in similar terms, using `the terrestrial
globe' in place of `the Earth':
\begin{quote}
[I]t would suffice here to explain the infinite through the
incomparable.  That is, to think of quantities incomparably greater or
smaller than ours.  This would provide as many degrees of
incomparability as we may wish.  Since that which is incomparably much
smaller has no value whatever in relation to the calculation of values
which are incomparably greater than it.  It is in this sense that a
bit of magnetic matter which passes through glass is not comparable to
a grain of sand, or this grain of sand to the terrestrial globe, or
the globe to the firmament.  
\cite{Le02}
\end{quote}
In sum, these texts by Engels make it clear that, far from viewing
infinitesimals as chimeras, he viewed them as a fine accomplishment
for which ``Nature provides the models.''

On the other hand, in the correspondence between Marx and Engels in
1881--82, Engels recognized the innovative value of Marx's
contribution.

\section{Synthesis: Pascal--Leibniz characteristic triangle}
\label{s4}

Marx's correspondence with Engels reveals his evolving views on the
calculus.  He initially endorsed infinitesimals but later adopted a
more critical stance.  In the 1860s,%
\footnote{Yanovskaya dated this text at 1865/6.  
\cite{An25} dates this text at 1863 based on a 2010 edition of Marx
and Engels.}
Marx sent to Engels the following explanation of the determination of
the tangent line to a parabola:
\begin{quote}
If I drop a perpendicular~$np$ to the axis, then~$p$ must be
infinitely close to~$P$, and~$np$ an infinitely close parallel line to
$mP$.  Now drop an infinitely small perpendicular~$mR$ to~$np$.  If
you now take the abscissa~$AP$ for~$x$, and ordinate~$mP$ for~$y$,
then~$np=mP$ (or~$Rp$), increased by an infinitely small increment
$[nR]$ or~$[nR]=dy$ (differential of~$y$), and~$mR= (Pp) =dx$.  Since
the part~$mn$ of the tangent is infinitely small, it coincides with
the corresponding part of the curve itself.  \cite[252; translation
  ours]{Ma68}
\end{quote}
\begin{figure}
\begin{center}
\includegraphics
[scale=0.5,
trim = 0 0 0 0,clip]
{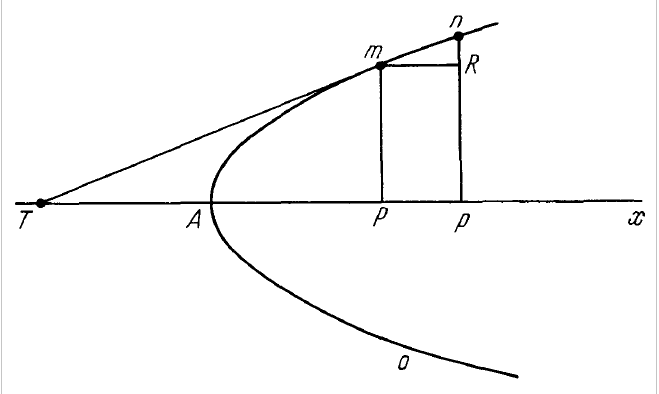}    
\end{center}
\caption{Determination of tangent to parabola}
\label{parabola}
\end{figure}
See Figure \ref{parabola}.  Marx goes on to assert that
triangles~$mnR$ and~$mTP$ are similar.  Here~$mnR$ is the
infinitesimal \emph{characteristic triangle} Leibniz learned from
Pascal.  For the purpose of this calculation, the arc $mn$ is assumed
to be a straight-line segment.  Infinitesimals are mentioned six times
in this passage.  Here Marx takes for granted the usefulness of
infinitesimals in explaining geometric properties of conic sections.
Infinitesimals are neither avoided nor described as
chimerical/mystical.

This important text was postponed until Part II of Yanovskaya's
edition of Marx's mathematical works \cite{Ma68}, and was not included
in the English version \cite{Ma83} (even though three letters between
Marx and Engels dating from 1881/2 are included prominently in the
\emph{Mathematical Manuscripts of Karl Marx} \cite[xxvii--xxx]{Ma83}
even \emph{before} Marx's mathematical manuscripts themselves).

Marx learned this infinitesimal-based approach from \cite{Sa74}.  The
title of this section is, of course, tongue-in-cheek.  Marx went on to
study \cite{Bo27} who (for the most part) avoided infinitesimals.
Marx subsequently adopted a progressively more hostile attitude toward
infinitesimals, as illustrated in Section~\ref{s1}.  Influenced by
Boucharlat and others, he eventually adopted an approach analyzed in
Section~\ref{s3}.

\section{Marx on preliminary derivatives and~$0/0$}
\label{s3}

Let us analyze Marx's comments on the definition of the derivative.
But first, a word on crises of capitalism.

\subsection{Calculus and periodic crises}

As part of his commitment to social activism, Marx may have been
interested in the calculus because of his interest in proving (perhaps
by analyzing a suitable differential equation?) that `unplanned'
capitalist economy inevitably leads to periodic economic `crises'.%
\footnote{See e.g., \cite[219--220]{Ma83} for some comments on such
`periodic crises'.}
He was in correspondence with Samuel Moore (1838--1911) about the
matter, and wrote to Engels as follows: ``I have repeatedly attempted,
for the analysis of crises, to compute these `ups and downs' [in
  prices] as fictional curves, and I thought (and even now I still
think this possible with sufficient empirical material) to infer
mathematically from this an important law of crises.  Moore, as I
already said, considers the problem rather impractical, and I have
decided for the time being to give it up'' (Marx to Engels
\cite[x]{Ma83}).
Moore apparently felt that it would be impossible to prove such a
thing because there are too many variables involved.  Marx persisted
and even discussed this with Engels repeatedly (both in 1873 and
1882).

Vogt claims that ``The chemist Carl Schorlemmer told [Marx] that it
might be possible with the aid of calculus – more specifically, with
differential calculus – to calculate when the next crisis would
come.''%
\footnote{See
\url{https://libcom.org/article/marx-and-mathematics-annette-vogt} }
However, she provides no primary source to substantiate her claim.

\subsection{Carchedi on the Marxian approach}
\label{s52}

As far as the calculus is concerned, according to Carchedi, Marx
proceeds as follows:
\begin{quote}
Given a certain function, such as~$y=f(x)$, Marx first lets~$x_0$
become~$x_1$.  Both~$x$ and~$y$ increase by \emph{finite}
quantities,~$\Delta x$ and~$\Delta y$ (so that the rules for ordinary
numbers can be applied here).  The ratio
\begin{equation}
\label{21}
\frac{\Delta y}{\Delta x}=\frac{f(x_1)-f(x_0)}{x_1-x_0}
\end{equation}
is what he calls the provisional or preliminary derivative.  Then, he
lets~$x_1$ return to~$x_0$ so that~$x_1-x_0 = 0$ and thus~$y_1-y_0=0$,
thus reducing this limit value to its absolute minimum quantity.  This
is called the definitive derivative,~$dy/dx$.%
\footnote{\label{f8b}
\cite[421]{Car}.  For further analysis of Marx's comments on the
calculus, see
%
\cite[112--113]{Ku91}.}
\end{quote}
Let us examine Marx's procedure in more detail.  In \cite[6]{Ma83} he
uses the example of a general quadratic function to introduce the
concept of a \emph{preliminary `derivative'} expressed in terms
of~$x_1-x=\Delta x$, corresponding to the difference quotient
\emph{before} the `higher-order terms' are discarded.  By page 7, he
introduces the term \emph{`derived' function} (i.e., the ordinary
derivative).  By page 8, he mentions ``functions of~$x$ which can only
be represented by infinite series, which [is] all too often the
case.''  Thus he implicitly acknowledges that the procedure he
outlined works only for analytic functions.

A special role in Marx's analysis is played by the symbol~$\frac00$.%
\footnote{Moigno introduced the term \emph{indeterminate form} for
$\frac{0}{0}$ familiar from the calculus.  Since Marx was in
possession of Moigno's textbook (see note~\ref{f2}), Marx's use of
$\frac00$ may have been inspired by Moigno.  To our knowledge,
Moigno's possible influence on Marx has not been studied in detail,
though there is some discussion in \cite{Al19}.}
Concerning the function~$f(x)=x^m$, he writes:
``{\ldots}~setting~$x_1-x=0$, or~$x_1=x$, turned~$\frac{y_1-y}{x_1-x}$
into~$\frac00$, and we substituted~$\frac{dy}{dx}$ in the place
[of~$\frac00$] to show what the origin of this~$\frac00$ is, i.e.,
which ratio of actual differences -- in the above case
$\frac{y_1-y}{x_1-x}$, finally turns into~$\frac00$.  This is all the
more justified since as a result we obtain
\[
\frac00=mx^{m-1}=f'(x),
\]
etc.''  \cite[474; translation ours]{Ma68}.  Such an analysis would
not get an A on a freshman calculus exam.

\subsection{Yanovskaya on Euler and Marx}

Editor Yanovskaya claims in \cite[485]{Ma68} a similarity between
Marx's and Euler's view of differentials as `quantitatively zero'.
However, as she acknowledges in \cite[578]{Ma68}, Euler worked with a
pair of `equalities', arithmetic and geometric.  Euler's arithmetic
equality between terms~$T$ and~$T'$ works as a relation of infinite
proximity, meaning that the difference~$T-T'$ is infinitesimal;
whereas the geometric equality entails that the ratio~$\frac{T}{T'}$
is infinitely close to~$1$.  For further details see~\cite{17b}.
There is no trace of such a distinction in Marx.

Marx's distinction between the \emph{preliminary derivative} and the
\emph{derived function} (see Section~\ref{s52}) is useful only in
explaining the differentiation of analytic functions, but not more
general functions.  To give an example of a type already mentioned by
Cauchy, consider the function~$f$ which vanishes at~$0$ and is
otherwise given by~$f(x)=e^{-\frac{1}{x^2}}_{\phantom{I}}$.  Since~$f$
tends to zero faster than any polynomial as~$x\to0$, we have
$f'(0)=0$, and similarly for all the derivatives at the origin (though
this is not something that the techniques Marx uses would enable him
to determine).  Thus the Taylor series in $x$ identically vanishes.

The `preliminary derivative' at the origin in this case would be
\[
\frac{e^{-\frac{1}{(\Delta x)^2}}}{\Delta x_{\phantom I}}.
\]
Since in this case one cannot ``cancel out'' the~$\Delta x$ in the
numerator and the denominator, Marx's distinction between the
preliminary derivative and the derived function sheds no light on the
``terms that need to be discarded'' in this case.  Indeed, ``Readers
will find these mathematical manuscripts to be a conceptual and
tropical desert worthy of only a cursory glance'' \cite[234]{Ol84}.

\section{Fermat, adequality, and elusive increments}
\label{s6}

How original is Marx's procedure as analyzed in Section~\ref{s3}?  It
is reminiscent of the technique of adequality as developed by Pierre
de Fermat two and a half centuries earlier; see \cite{Fer}.  We don't
know whether Marx was familiar with Fermat's writings, or whether
there was a direct influence by Fermat on Marx.  Rather, these ideas
were entirely conventional and customary for their time (as already
pointed out by Gumbel and condemned by the `Marx--Engels Institute'),
and were so to speak `in the air'.  The point is that Fermat's work
was even earlier than that of the inventors of the calculus, who were
themselves influenced by Fermat and mentioned such influences in their
writings.

Fermat proceeds as follows in applying his method of adequality to the
solution of the problem of determining the maximum of an expression
which we will denote~$f$ (of course Fermat did not use functional
notation).  Compare the values~$f(A)$ and~$f(E)$ by \emph{adequality}:
\[
f(A)\adequal f(E).%
\footnote{\label{f9}The notation~``$\adequal$'' is Leibniz's, not
Fermat's; see
%
\cite{KS1}.  Fermat himself used the term \emph{adequality} either in
Latin or French, according to the language of the text.  We
use~$\adequal$ to simplify the formula and remind the reader that
Fermat did \emph{not} use any notation for the equality sign to signal
his relation of adequality.  Rather, he used semi-symbolic notation
inspired by Vieta.}
\]
Express~$f(E)$ in terms of~$f(A)$ and~$A-E$, cancel out common terms,
and divide by~$A-E$.  In modern terms this corresponds, roughly, to
forming the quotient~$\frac{f(A)-f(E)}{A-E}$; cf.~Marx's procedure
\eqref{21} as reported by Carchedi.

Finally, simplify the resulting expression by clearing denominators,
and set~$A=E$.  This produces the required extremum.  We see that
Fermat's procedure fits well with Marx's procedure as outlined above.

In other texts, Fermat denotes the \emph{increment} by~$E$.  He then
compares~$f(A+E)$ and~$f(A)$ by adequality
\[
f(A+E)\adequal f(A),
\]
cancels out common terms, divides by~$E$ (or, as he says, a higher
power of~$E$) and discards higher-order terms.%
\footnote{Additional details on Fermat's method can be found in
%
\cite{An};
%
\cite{Giu};
%
\cite{KSS};
%
\cite{De17};
%
\cite{18d}.
}
Fermat's increment~$E$ played the role of the Leibnizian differential
$dx$ (which was the focus of Marx's analysis) in the procedures of the
infinitesimal calculus.

\subsection{Fermat's algorithm}

Fermat's text \emph{Methodus ad Disquirendam Maximam et Minimam} opens
with a summary of an \emph{algorithm} for finding the maximum or
minimum value of an algebraic expression in a variable~$A$.%
\footnote{Following the conventions established by Vieta, Fermat uses
capital letters of vowels A, E, I, O, U for variables, and capital
letters of consonants for constants.}
The algorithm can be broken up into seven steps as follows:
\begin{enumerate}
\item
Introduce an auxiliary symbol~$E$, and form~$f(A+E)$.
\item
Set \emph{adequal} the two expressions:~$f(A+E) \adequal f(A)$.
\item
\label{cancel}
Cancel the common terms on the two sides of the adequality.  The
remaining terms all contain a factor of~$E$.
\item
Divide by~$E$ (see also next step).
\item
In a parenthetical comment, Fermat adds: ``or by the highest common
factor of~$E$".%
\footnote{For a discussion of this step, see
%
\cite[Section~3.5]{20e} following \cite{Giu}.}
\item
\label{among}
Among the remaining terms, suppress all terms which still contain a
factor of~$E$.
\item
Solving the resulting equation for~$A$ yields the extremum of~$f$.
\end{enumerate}
Given that the algorithmic nature of the procedures involved was
already evident to the seventeenth century masters, it is surprising
to find recent authors attributing to Marx an alleged innovation of
viewing the procedures of the calculus as algorithms, as when Ricci
claims that ``[Marx's] notion is strikingly similar to the modern
concept of algorithm, making Marx a precursor of the modern
computational mathematics'' \cite[Abstract]{Ri18}.  Similar
algorithmic claims appear in \cite[40]{Ma21}.  Chime in Brenner and
Igamberdiev: ``These views of Marx anticipate the modern concept of
the algorithm, making him a precursor of modern computational
mathematics and logic.  {\ldots}~The operational nature of mathematics
was developed in the twentieth Century in the logical approaches of
intuitionism and constructivism'' \cite[155]{Br21}.  Attempts to
connect Marx to intuitionism are surprising as the fundamental feature
of intuitionism is the rejection of the Law of Excluded Middle.  There
is no explicit source in Marx for such a development.

Marx considered the approximation of the derivative~$\frac{dy}{dx}$ by
the ratio~$\frac{\Delta y}{\Delta x}$ to be conceptually flawed, since
for finite non-infinitesimal~$\Delta x$, the ratio will only
approximate the derivative.  However, for infinitesimal~$E$,
discarding higher-order terms in~$\frac{f(A+E)-f(A)}{E}$ results
precisely in the value of the derivative at~$A$ (and not merely an
approximation).  Fermat (who did not possess the concept of
derivative) was very cautious in his discussions of the nature of~$E$
(and certainly did not refer to it as an infinitesimal); the possible
reasons (some of them theological) are analyzed in \cite{18d}
and~\cite{23g}.  In modern infinitesimal analysis, the procedure of
``discarding higher-order terms" is formalized in terms of taking the
standard part, as discussed for example in the textbook~\cite{Ke}.

\section{Kolman: ``Marx really differentiates''}
\label{s7}

The source of attributing prodigious mathematical talent to Marx can
be pinpointed to a specific Soviet Apparatchik.  This is Ernst
(``Arnosht'') Kolman (1892--1979) in the year 1931.  Writes Seneta:
\begin{quote}
According to his biography [Kolman, 1982b] he did not return to Moscow
and academic activity {\ldots}~till March, 1931, at the
Marx--Engels--Lenin Institute (MELI), where he was put in charge of
the Marx office.  The recently arrested director of MELI, Riazanov,
was thought to have suppressed publication of many of Marx's
manuscripts, in particular Marx's mathematical manuscripts, since
Riazanov was not convinced of their academic value, on the advice, as
Kolman p.\,172 says, of the ``mediocre German mathematician Gumbel%
\footnote{Brenner notes: ``Gumbel published a summary of the papers in
`O matematitscheskich rukopisach K. Marxa', \emph{Letopisi Marxisma} 3
(1927), but according to Vogt (20--23), publication of the documents
themselves never followed, presumably because they were inconsistent
with the mythology of Marx as genius'' \cite[182, n.\;99]{Br93}
(we preserved Brenner's transliteration of the Russian).}
{\ldots}\;who did not grasp their methodological significance.''
Something similar, says Kolman [1982b, 172] occurred with Engels'
\emph{Dialectic of Nature}.  Einstein was asked to comment on the
significance of this but ``{\ldots}~did not grasp the enormous
philosophical significance {\ldots}~since Einstein generally did not
understand dialectics.'' \cite[360]{Se04}
\end{quote}
With regard to \emph{Dialectic[s] of Nature}, we will take Einstein's
word for it, \emph{pace} Kolman.%
\footnote{Kolman apparently held a dim view of Einstein's physics, as
well.  Thus, ``Kolman (1939) believed that velocities can exceed
300,000 km/sec.  The contrary statement, he declared, went against
dialectical materialism'' \cite{Sh03}.}
In Kolman's, the disgraced Riazanov (see further in Section~\ref{s81})
did not sufficiently appreciate Marx's mathematical talents, but
Kolman did.

\subsection{Dialectical Materialism 101}
\label{s71b}

A reader unfamiliar with this circle of ideas may find the following
summary helpful as a zeroth approximation.

In this scheme of things, \emph{idealism} is mystical and bad, whereas
\emph{materialism} and \emph{dialectics} are good.  Hegel was
dialectical but idealist, D\"uhring was materialist but not
dialectical.  Marx, Engels, and Lenin were both (materialist
\emph{and} dialectical).

To strip away the mystical veil of the calculus, a formal-logical
foundation is insufficient; dialectical logic resolves the paradox of
infinitesimals by viewing them not as static, mystical entities but as
a representation of a process of becoming and change, a unity of
opposites between being and non-being.  Marx saw the foundational
problem of the calculus as a `touchstone' for applying dialectical
materialism to mathematics; see e.g., \cite[316]{Ke77}.  Modern set
theory (as practiced by Luzin) was a foe of Marxist epistemology
\cite[67]{Vu00}.  The formalist orientation in the foundations of
mathematics, which favored the axiomatization of mathematics, was an
idealistic delusion \cite[34]{Ya50}, \cite[69]{Vu00} or idealistic
aberration \cite[30]{Vu02}.

Vucinich explored the Soviet endeavor to create a `Marxist
mathematics', a concept that sought to align mathematical thought with
the tenets of dialectical materialism, including the Soviet critique
of certain areas of modern mathematics as `idealist' or `bourgeois',
and the promotion of a mathematics that was seen as more aligned with
materialist philosophy.  Logicism, formalism, and intuitionism are the
most influential orientations in the foundations of mathematics that
must be rejected, being manifestations of philosophical idealism
\cite[108]{Vu99}.

How long did such a `delusion' last?  Vucinich noted that by 1991,
``the [Soviet] journal \emph{Priroda} (\emph{Nature}) published a
Russian translation of an article by Jean van Heijenoort%
\footnote{The name is misspelled in \cite{Vu02}.}
who assembled enough arguments to conclude that Engels knew little
about contemporary mathematics, that his `materialism' was barely more
than `crude empiricism', and that his dialectics was a `degenerate and
oversimplified' notion of Hegel’s dialectics'' \cite[36]{Vu02}.  The
original essay appeared as \cite{He85} but most of it dates from 1948.

\subsection
{Rhetoric of \emph{dialectics} and \emph{negation of the negation}}
\label{s71}

In 1968, Kolman opined that ``Marx, like Hegel, considered all efforts
to provide a purely formal-logical foundation for analysis hopeless
{\ldots}~He set himself the task of providing a foundation for
analysis \emph{dialectically}, relying on the unity of the historical
and logical aspects'' (Kolman as translated in \cite[227]{Ma83};
emphasis added).  Apparently, `dialectics' is capable of dressing up
the bug (of the absence of a `formal-logical foundation') as a feature
(see Section~\ref{s71b}).  Naturally, therefore, ``in the development
of mathematics the formal-logical moments of understanding are
shouldered aside by the dialectical moments'' (Kolman and Yanovskaya
as translated in \cite[241]{Ma83}; in such literature, the term
\emph{moment} is often used in the sense of \emph{aspect}).
Accordingly, dialectician Hegel is credited in the following terms:
``[Hegel] gave an essentially materialistic definition of mathematics
which \emph{smashes apart} the framework of the bourgeois
world-outlook with its characteristic quantity fetishism.''%
\footnote{\cite{Ko31}; emphasis added.  See \cite[236--237]{Ma83} and
\url{https://www.marxists.org/reference/subject/philosophy/works/ru/kolman.htm}.
The importance of materialism was reiterated by Struik 17 years later:
``Marx' position was that of the materialist, who insists that
significant mathematics must reflect operations in the real world''
\cite[193]{St48}.}
Many mathematicians' careers were indeed shouldered aside and smashed
as a result of social activism by the likes of Kolman and his coauthor
Yanovskaya; see \cite[131]{Ba01} and \cite{Vu00}.

Ricci holds a similarly high opinion of dialectics as a tool to
account for~$\frac{dy}{dx}$:
\begin{quote}
From a methodological point of view, the definition of the increment
in terms of negation of the original variable sets in motion a
\emph{dialectical development} that, through the subsequent negation
of the negation, will lead to the final result. {\ldots}~The
increments becoming zero, which in \emph{dialectical terms} represents
the \emph{negation of the negation}, only takes effect on the
left-hand side of the equation, which is reduced to the
expression~$0/0$, leaving the right-hand side unchanged.  However,
this time the ratio~$0/0$ should no longer be dreadful because it does
not denote an arithmetic operation, but is a purely symbolic operator,
which as such can be replaced by the differential ratio~$dy/dx$,
without giving rise to any logical contradiction.  \cite[224--225;
  emphasis added]{Ri18}
\end{quote}
While rich in metaphor, such flourishes of the rhetoric of
\emph{dialectics} and \emph{negation of the negation} would not get an
A on a calculus test, either.

\subsection{A. D. Aleksandrov on dialectics}

The problem with attributing to dialectical materialism the ability to
solve specific mathematical problems (such as the definition of the
derivative) was pinpointed by a well-known advocate of dialectical
materialism, A.\;D.\;Aleksandrov, in the following terms:
\begin{enumerate}
\item
Dialectical materialism, needless to say, does not offer methods
for solving specific problems in mathematical science, but
\item
it indicates true reference points for searches for scientific truth
and arms one with methods for elucidating the true import of theories
and the content of scientific concepts.%
\footnote{Aleksandrov as quoted in
%
\cite[263--264]{Gr87}; numerals (1) and (2) added.}
\end{enumerate}
Regardless of one's opinion of Aleksandrov's item (2), his item (1)
needs to be taken to heart by Carchedi, Ricci and others.  The
evolution of Aleksandrov's views is analyzed in \cite[26]{Vu02}.

\subsection{New Foundation for the calculus}

Kolman appears to have initiated the trend of attaching superlative
mathematical significance to these manuscripts.  In a 1932 lecture, he
claimed that we owe Marx no less than `A new foundation of the
Differential Calculus'
%
\cite{Ko32}.  A further insight is that
``[Marx] \emph{really differentiates}, thanks to which differential
\emph{symbols} appear, while Lagrange applies differentiation to the
algebraic binomial expansion.''%
\footnote{Kolman \cite[227]{Ma83}; emphasis on ``really
differentiates'' added.}

Kolman's approach reflected his attitude of worship toward the
ideologues of communism.  In fact, the manuscripts are insignificant
mathematically (by post-seventeenth century standards), but they do
reflect Marx's interest in seeking novel mathematical approaches to
social and economic problems.

\section{Editorial history from Gumbel to Yanovskaya}
\label{s8}

%
\cite{Vo19} traces the history of attempts to edit Marx's mathematical
texts and the fate of the successive editors.  The first editor was
Emil J. Gumbel (1891--1966), who prepared a commentary on a pair of
mathematical manuscripts by Marx already in 1927, and tried repeatedly
to get them published in Moscow.%
\footnote{In the 1983 edition, we find the following comment: ``To
work on them the Institute commissioned the German mathematician
E. Gumbel, whom R. Matejka and R. S. Bogdan helped to decipher the
extremely difficult text'' \cite[225]{Ma83}.  Actually R. Matejka
and R. Bogdan (maiden name) is one and the same person.}
This was delayed year after year, awaiting final approval which never
came.

\subsection{Purge of Riazanov}
\label{s81}

The editors in Moscow in charge of Gumbel's submission were David
Riazanov (1870--1938) and his subordinates, including Ernst Cz\'obel.
Riazanov was the director of the Marx--Engels Institute and the
initiator and head of the project of publishing the collected papers
of Marx and Engels, who handled the negotiations with the Archive of
the German SPD and arranged for a permission enabling Moscow to obtain
photocopies of the papers of Marx and Engels so as to start with the
publishing project.

Unfortunately for both Gumbel and Riazanov, Riazanov had a falling-out
with the ruling elite in 1931, was arrested, and stripped of all his
honors.  He eventually stood a sham trial and was executed.  In
retrospect, Riazanov's purge was inevitable, given the animosity
between him and Stalin dating back to 1921; see \cite[6]{Ro08}.

\subsection{Conscientious work}

The person who helped Gumbel decipher some of Marx's manuscripts was
R\'ozsa Matejka or Mathejka (n\'ee Bogdan).


In 1933 three manuscripts were finally published by a group of three
editors including Yanovskaya.  Yanovskaya acknowledged that already
``by 1927, a significant part of the material had been deciphered and
\emph{quite conscientiously} at that.''%
\footnote{Emphasis added.  Transliterated original: ``K 1927 godu
znachitel'naya chas't materiala byla rasshifrovana pritom vpolne
dobrosovestno''
%
\cite[75, note\;2]{Ya33}.  Note 2 goes on to castigate the unnamed
decipherers for insufficient appreciation of Marx's manuscripts, and
even mentions Gumbel's 1927 publication \cite{Gu27} without however
mentioning his name.}
Note the passive voice.  Neither Gumbel nor Matejka/Bogdan are
mentioned here by name (though Bogdan is mentioned forty pages later
in a footnote on the last page of the article).  Exploitation of other
scholars' conscientious work without mentioning their names amounts to
plagiarism.

Nearly a century later, Antonova claims the following:
\begin{quote}
{\ldots}~apart from the publication of E. J. Gumbel's article `On the
Mathematical Manuscripts of K. Marx (Communication)' (Letopisi
Marksizma [Chronicles of Marxism].  Notes of the Institute of K. Marx
and F. Engels.  Moscow--Leningrad: State Publishing House,
1927. V. 3), 56--60 and A.\,Vogt's article on the preliminary stage of
work on systematizing the manuscripts and preliminary plans for their
publication, nothing is known about this work in principle for the
sole reason that the work did not go beyond preliminary plans.
\cite[10]{An25}
\end{quote}
Antonova's claim that Gumbel's editorial work did not go beyond
preliminary plans is in error.  In fact, \cite{Vo19} shows that the
work reached the galley proof stage but was never published for
reasons beyond the control of editor Gumbel, as explained above.
Antonova's claim that nothing is known about this work is similarly
inaccurate.  Indeed, as explained above, Yanovskaya commented in 1933
concerning Gumbel and Matejka's \emph{conscientious work} of
deciphering several manuscripts (without mentioning either name).
Yanovskaya's \emph{conscientious work} comment is not analyzed by
Antonova in her discussion of the 1933 edition in \cite[23--28]{An25}.

In a 2024 article, Antonova quotes at length from a 1927 letter sent
to Gumbel from Moscow, and concludes as follows: ``Already from this
letter of July 21, 1927 it becomes clear that D.B. Ryazanov was not
satisfied with Gumbel’s interpretation of mathematical manuscripts,
and no longer envisaged Gumbel's participation in the work on the
publication of Marx's mathematical manuscripts'' \cite[70]{An24}.
However, in 1930, Riazanov wrote to Gumbel in the following terms: ``I
can, however, reassure you that no one will obstruct the publication
or handling of these studies of Marx {\ldots}.  This work will
definitely get published in the next - third - German volume of the
Archiv.  This volume will certainly appear within this working year,
by the summer of 1931.  I can assure you that the institute, which got
the editorial work {\ldots} by you, has no less interest in the
publication than you yourself have'' \cite[28]{Vo19}.  It is therefore
surprising to find Antonova claiming that ``The reason that
E.J. Gumbel's article in the `Archive of K. Marx and F. Engels' was
not published was that D.B. Ryazanov's consent for its publication was
not obtained'' \cite[70]{An24}.

Four months later, Riazanov was arrested, and the Marx--Engels
Institute purged of his supporters; see \cite{Ro08}.  Were it not for
stalinist purges, Gumbel's article would surely have appeared, as
promised by Riazanov.

Some later editors, including Yanovskaya and Antonova, disagreed with
Gumbel's approach to the mathematical manuscripts of Marx.%
\footnote{In a 1953 internal report \cite{Ya53}, Yanovskaya went as
far as describing Gumbel's work as `a libel on Marx's manuscripts'
(Yanovskaya as quoted in \cite[84, note 132]{An24}).}
However, such disagreement does not justify suppressing his name from
the history of the editorship of some of the manuscripts.

\subsection{Kolman's reappearance}

The good stalinist of the 1930s, Kolman was himself canceled in 1968
for expressing displeasure with Soviet actions in Czechoslovakia.
Even though he had been a driving force behind the manuscript project
from the start, Yanovskaya's coauthor Kolman was not even mentioned in
the 1968 edition \cite{Ma68}, as noted by Vogt (nor is Gumbel
mentioned).  Kolman reappears in the 1983 English edition.

Yanovskaya died in 1966, two years before the book was published,
allegedly edited by her.
The 1983 English edition~\cite{Ma83} cautiously lists no editor at
all.

The manuscripts show that Marx respected and appreciated mathematics,
but possess little scientific value.  The technique of differentiation
as he describes it works only for polynomials (modulo jettisoning
dialectical material like~$\frac00=\frac{dy}{dx}$); even if one
credits him with wishing to extend it to power series, it is still
limited to analytic functions, similarly to the approach in early
Lagrange.  It is typical of the epoch that embellishing the
\emph{cul-de-sac} of Marx's approach to the calculus is linked to the
personality of Kolman, known for his infamous role in the Luzin
affair, including sinister accusations of `Idealism'.%
\footnote{See \cite[33--34]{Vu02}, \cite{GK}; \cite{Ku12};
\cite{Go17}; \cite{Ne21}.}
%
%

\section{Whatever became of Leibnizian infinitesimals?}
\label{s9}

We saw in Section~\ref{two} that the dialecticians in the vanguard of
the proletariat differed as to the utility of Leibniz-style
infinitesimals.  Ironically, just as Engels was heaping praise upon
them, the ``great triumvirate'' \cite[298]{Boy} of Cantor, Dedekind,
and Weierstrass sought skillfully, though not entirely successfully,
to maneuvre them out of existence.%
\footnote{See
%
\cite{Eh06},
%
\cite{KS2}, \cite{KS1},
%
\cite{24b}.}

In the second half of the twentieth century, building on earlier work
by \cite{He48} and \cite{Lo55}, Robinson in \cite{Ro61} and
\cite{Ro66} developed a formalisation of infinitesimals usable in
analysis and other branches of mathematics.  Robinson's work was
anchored in a formal-logical foundation and modern set theory (both
deprecated by dialectical materialism), in a remarkable validation of
neo-Kantian Hermann Cohen's intuitions about infinitesimals.

A key insight is that the Leibnizian heuristic principle, ``the rules
of the finite succeed in the infinite and vice versa,'' admits a
formalisation as the Transfer Principle of nonstandard analysis (this
insight is quite at variance with the claim in \cite[113]{Ma83} that
the algebra of the calculus of Newton and Leibniz is ``as different
from the usual algebra as Heaven is wide'').  In this way, ``Robinson
succeeded in showing the reasonableness of `redrawing' the early
history of the calculus to reinstate past views that, cast in the
light of nonstandard analysis, could be seen more clearly''
\cite[327]{Da21}.  Robinson's philosophical position of Formalism is
clarified in \cite{26c}.

Following in the footsteps of Robinson's pioneering work of the 1960s,
\cite{Hr78} and \cite{Ne} developed axiomatic approaches to
nonstandard analysis, exhibiting infinitesimals within the real
line~$\R$ itself.  While modern set-theoretic \emph{foundations} are
alien to the conceptual framework possessed by Leibniz, the
\emph{procedures} of the axiomatic approach to nonstandard analysis
are arguably closer to the spirit of the Leibnizian calculus than the
procedures of the model-theoretic approach.

\section{Acknowledgments} We are grateful to historian of mathematics
Annette Vogt for helpful comments.  The influence of Hilton Kramer
(1928--2012), who pursued and exposed the mendacity of the marxist
`long march' in academia and the arts with mathematical rigor, is
obvious throughout.

\end{document}